\documentclass[a4paper,11.5pt,reqno]{amsart}

\usepackage{graphicx}
\usepackage[utf8]{inputenc}
\usepackage{amsmath,amssymb,amsthm}
\usepackage{enumitem}
\usepackage[dvipsnames]{xcolor}
\usepackage{hyperref}
\usepackage{caption}

\usepackage{lineno}
\modulolinenumbers[5]
\numberwithin{equation}{section}

\hypersetup{colorlinks = true, linkcolor = blue, urlcolor = blue, citecolor = blue}

\renewcommand{\theta}{\vartheta}

\begin{document}

\title[Cubic equations with 2 Roots in an interval]
      {Cubic equations with 2 Roots in \\ the interval $[-1, 1]$}

\author[H. Ruhland]{Helmut Ruhland}
\address{Santa F\'{e}, La Habana, Cuba}
\email{helmut.ruhland50@web.de}



\begin{abstract}
The conditions for cubic equations, to have 3 real roots and 2 of the roots lie
in the closed interval $[-1, 1]$ are given. These conditions are visualized. This question arises
in physics in e.g. the theory of tops.
\end{abstract}

\date{\today}

\maketitle

\section{Introduction}

Cubic polynomials are ubiquitous in physics. I cite here just some examples from the introduction in \cite{Prod}:   \textit{"\dots The applications of cubic and quartic equations in all branches of science are vast. \dots There are well over 200 real gas equations, many of which are also cubic. The elastic waves propagating on the surface of solids, the so-called Rayleigh waves \dots The Hodgkin-Huxley model in mathematical neuroscience encounters a quartic \dots In general relativity, through the d'Inverno and Russel-Clark algorithms, the Petrov classification of the Weyl conformal curvature \dots"} \\

The question, if a cubic polynomial has only real roots, can be decided using the discriminant condition.
This question can be extended to the question:
When has a cubic polynomial $0 \dots 3$ roots in a given open or closed interval? \\

This question arises e.g. in physics in the theory of tops, i.e. rigid bodies that move under the influence of  gravity around a fixed point with 3 degrees of freedom (the 3 Euler angles). Here in the case of nutation the upper and lower limits for $\cos (\theta)$ are given by $2$ real roots of a cubic in the closed interval $[-1, 1]$. The $3^{rd}$ real root lies outside this interval, see appendix \ref{ExPhysics}. \\

Description of the problem treated in this article: \\

\textit{Determine the conditions under which exactly 2 roots of a monic cubic polynomial
${x^3 + a x^2 + b x + c}$ lie in the closed interval $[-1, 1]$.} \\

\newpage

\section{The conditions for a cubic polynomial with $2$ roots in $[-1, 1]$ \label{Cond_one_look}}

Let $P = x^3 + a x^2 + b x + c$ be a monic cubic polynomial. The discriminant is defined as:
\begin{equation}
   D_3 = - 27 \,  c^2 + (18 \, a b - 4 \, a^3) \, c + a^2 b^2 - 4 \, b^3
   \label{def_D_3}
\end{equation}
Assume the condition $D_3 \ge 0$ to get $3$ real roots is fulfilled. \\
 
Calculate these $5$ quantities: 
\begin{equation}
   A = a + b + c + 1	\quad B = a - b + c - 1 \quad
   A_T = 4 \, (c + 1) \quad B_T = 4 \, (c - 1) \quad  E = (a - c) \, c - b + 1
   \nonumber
\end{equation}

Distinguish 3 cases depending on $c$: \\

{ \bf 1. $\mathbf {c < 0}$ } \\

Replace $a \rightarrow - a, \enspace c \rightarrow - c$, this is the map $M : x  \rightarrow - x$. Now $c > 0$ and we get one of the following cases. \\

{ \bf 2. $\mathbf {0 \le c \le 1}$ } \\
\begin{equation}
   \begin{split}
	 &( \; A < 0 \quad \enspace \;  and \enspace \;  B \le 0 \; )  \quad or \quad
        ( \; A \ge 0 \enspace  \; and \enspace \; B > 0   \enspace \, ) \quad  or  \\
	 &( \; A > A_T \enspace \;  and \enspace\;  B = 0   \; \, )  \quad or \quad
        ( \; A = 0   \enspace \; and \; B < B_T \; )
   \nonumber
   \end{split}
\end{equation}

{ \bf 3. $\mathbf {c > 1}$ } \\
\begin{equation}
	( \; A \le 0 \enspace \; and \enspace \;  B \le 0 \; )  \quad or \quad
   ( \; A \ge 0  \enspace \; and \enspace \; B \ge 0   \; \enspace \; and \enspace \; E \ge 0 \; )
   \nonumber
\end{equation}

\section{The discriminant surface $D_3 = 0$}

\subsection{The $2$ components of the discriminant surface}

In the $a-b$ plane the discriminant surface consists of a parabola for $c = 0$, it consists of $2$ components for $c \ne 0$: \\
For $c > 0$ a smooth component at the left of $P_b$ and below of $P_a$, the $2$ parabolas in the following figure \ref{comp2}. The parabolas intersect at $(0, 0)$ and have  perpendicular axes. The second component, smooth with the exception of a cusp is located inside the $2$ parabolas. All cusps lie on the parabola $P_C$:
\begin{equation}
   a_C = 3 \, c^{1/3}  \quad  b_C = 3 \, c^{2/3}  \qquad  P_C : a^2 - 3 \, b
\end{equation}

The $2$ parabolas are defined by:
\begin{equation}
   P_a = b^2 - 4 \, a c  \qquad  P_b = a^2 - 4 \, b
\end{equation}
The $2$ components approach in the limit $a, b \rightarrow + \infty$ to the parabolas. To see this e.g. for $P_b$ replace $b$ in the equation for $D_3$ by $a^2 / 4$. The two terms with $a^6$ cancel. The remaining terms are of size $O (a^4)$. So for $b \rightarrow + \infty$ the components 1 and 2 of $D_3$ approach to the parabola $P_b$. Component 1 from outside, component 2 from inside.
For $P_a$ replace $a$ in the equation for $D_3$ by $b^2 / (4 c)$ ...

\begin{center}
  \includegraphics[width=0.92\textwidth]{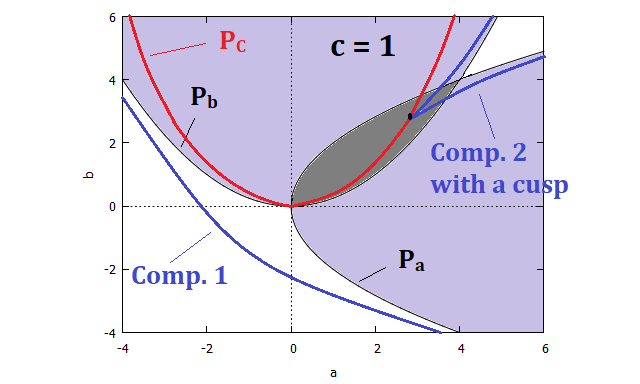}
  \captionof{figure}{$c = 1$, the 2 components of $D_3 = 0$, the 2 parabola $P_a, P_b$
                     The cusp is located in the dark grey shaded lens. The cusps for all $c$
                     lie on the red parabola $P_C$.}
  \label{comp2}
\end{center}

The parabola $P_C : a^2 - 3 \, b$ also shows up in the discriminant of the differentiated cubic
$D_2 = 4\, (a^2 - 3 b) = - 12 b^*$ with $b^*$ the second coefficient in the depressed cubic
$x^3 + b^* x + c^*$. For $3$ real roots besides $D_3 \ge 0$ this $D_2$ has to be $> 0$.
Figure \ref{comp2} visualizes that $D_3 \ge 0$ already implies $D_2 > 0$.

\subsection{The intersection of the planes $A$ and $B$ with the discriminant surface \label{inter}}

Define the following 2 planes:
\begin{equation}
   A = a + b + c + 1	\qquad B = a - b + c - 1
\end{equation}
The planes $A = 0$ and $B = 0$ represent polynomials with a root $+1$ or $-1$

\begin{equation}
   A_T = 4 \, (c + 1) \qquad B_T = 4 \, (c - 1) \qquad  A_{I \, 1/2} = 2 \, (c + 1 \pm 2 \, \sqrt{c})
\end{equation}
For $c \ge 0$: \\ 
The discriminant $D_3$ intersects with the plane $B = 0$ in $2$ lines $A = A_{I \, 1/2}$ (the subscript I means intersect) and $D_3$ is tangent to $B$ at the  line $A = A_T$ (the subscript T means tangent or touch). $D_3$ is tangent to the plane $A = 0$ at the line $B = B_T$  and doesn't intersect $A$ (the two $B_{I \, 1/2}$ are not real). \\   

The intersection of the 2 planes $A = 0$ and $B = 0$ in a line represent the polynomials $(x - 1) (x + 1) (x - c)$. \\

$A = 0$ and $B = B_T$ represent the polynomials $(x - 1)^2 (x + c)$, double roots because it's an intersection (tangent) with $D_3$. $B = 0$ and $A = A_T$ represent the polynomials $(x + 1)^2 (x + c)$. \\

$B = 0$ and $A = A_C$ represent the polynomials $(x + 1) (x - \sqrt{c})^2$, double roots because it's an intersection (though not tangent) with $D_3$. \\

$E$ is a ruled surface, for fixed $c$ a line. The lines $(A = 0, B = B_T)$ and $(B = 0, A = A_T)$ lie in
this surface $E$. Used in the condition \ref{Cond_cgt1}, figure \ref{c4} and \ref{c4_cusp}  for the case $c > 1$ to distinguish a different number of roots in the same quadrant:
\begin{equation}
   E = (a - c) \, c - b + 1
\end{equation}
The line $A_T, B_T$ is defined by $A \, / \, A_T + B \, / \, B_T - 1 = 0$. It follows $E = (A \, B_T + B \, A_T - A_T \, B_T) / 8$.

\section{The cubic polynomials with $2$ roots in the interval}

Colors in the following figures: \\
the discriminant in{ \color{blue} blue} (for $c = 0$ a parabola and a double line $b = 0$) \\
numbers in{ \color{red} red} show the number of roots in the interval \\
{ \color{LimeGreen} light green} shaded open regions inside the $2$ quadrants built by the lines $A$ and $B$ \\
{ \color{ForestGreen} dark green} open lines, the corresponding polynomial has $2$ roots in the interval \\
{ \color{ForestGreen} dark green} bullet, end point of a closed line with $2$ roots in the interval \\

\subsection{The case $c < 0$}

Replace $a \rightarrow - a, \enspace c \rightarrow - c$, this is the map $M : x  \rightarrow - x$. Now $c > 0$ and we get one of the following cases.

\subsection{The case $c = 0$}

Looking at the following figure \ref{c0}, we get this condition:
\begin{equation}
   \begin{split}
	 &( \; A < 0 \quad \enspace \;  and \enspace \;  B \le 0 \; )  \quad or \quad
        ( \; A \ge 0 \enspace  \; and \enspace \; B > 0   \enspace \, ) \quad  or  \\
	 &( \; A > A_T \enspace \;  and \enspace\;  B = 0   \; \, )  \quad or \quad
        ( \; A = 0   \enspace \; and \; B < B_T \; )
   \label{Cond_c01}
   \end{split}
\end{equation}

\begin{center}
  \includegraphics[width=0.92\textwidth]{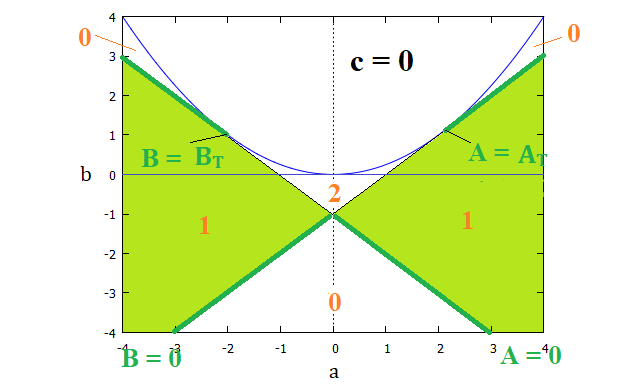}
  \captionof{figure}{$c = 0$, i.e. a root $0$ and roots of the quadratic polynomial $x^2 + a x + b$,
                     in red the number of roots of the quadratic in the interval $[-1, 1]$,
                     the green lines represent polynomials with 1 root in the interval
                     the blue curve is the parabola for the discriminant $D_2 = 0$.
}
  \label{c0}
\end{center}

\newpage

\subsection{The case $0 < c < 1$}

The condition is the same \ref{Cond_c01} as for the previous case, see figures \ref{c1_4} and \ref{c1_4_cusp}.
\begin{center}
  \includegraphics[width=0.92\textwidth]{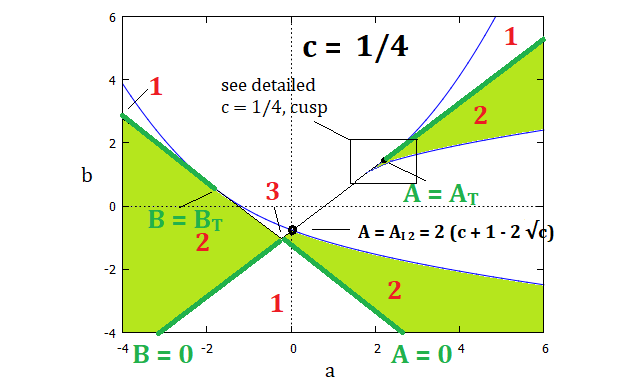}
  \captionof{figure}{$c = 1 / 4$, roots of the cubic polynomial $x^3 + a x^2 + b x + c$,
                     in red the number of roots in the interval $[-1, 1]$. The point $A = A_{I \, 2}$
                     is the intersection of the line $B$ with the discriminant $D_3$.}
  \label{c1_4}
\end{center}

\begin{center}
  \includegraphics[width=0.92\textwidth]{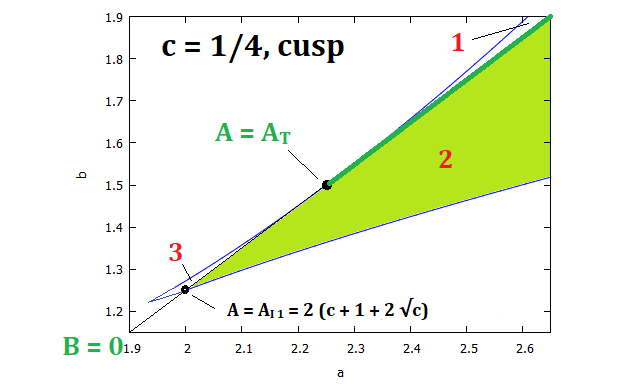} 
  \captionof{figure}{$c = 1 / 4$, the cusp, roots of the cubic polynomial $x^3 + a x^2 + b x + c$,
                     in red the number of roots in the interval $[-1, 1]$. The point $A = A_{I \, 1}$
                     is the intersection of the line $B$ with the discriminant $D_3$.}
  \label{c1_4_cusp}
\end{center}

\newpage

\subsection{The case $c = 1$}

The condition is the same \ref{Cond_c01} as for the two previous cases, see figures \ref{c1} and \ref{c1_cusp}.
\begin{center}
  \includegraphics[width=0.92\textwidth]{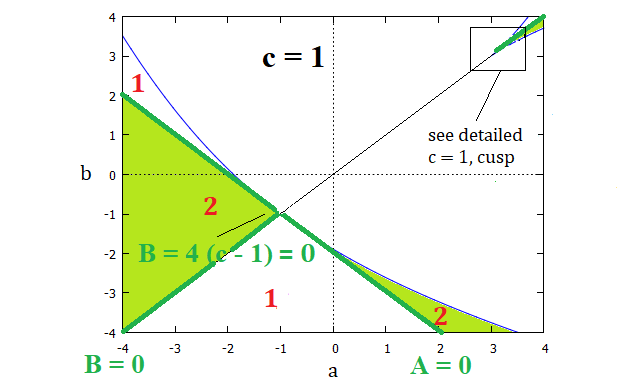} 
  \captionof{figure}{$c = 1$, roots of the cubic polynomial $x^3 + a x^2 + b x + c$,
                     in red the number of roots in the interval $[-1, 1]$.}
  \label{c1}
\end{center}

\begin{center}
  \includegraphics[width=0.92\textwidth]{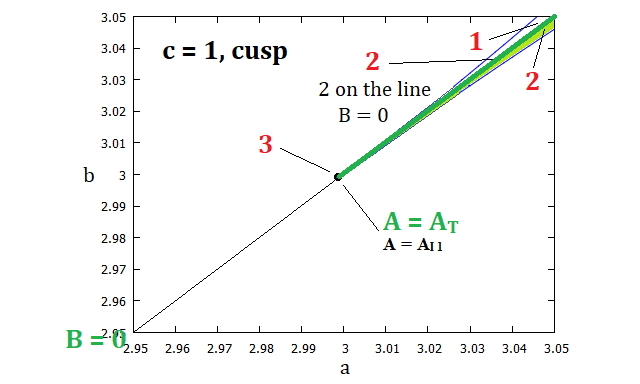} 
  \captionof{figure}{$c = 1$, the cusp, roots of the cubic polynomial $x^3 + a x^2 + b x + c$,
                     in red the number of roots in the interval $[-1, 1]$.}
  \label{c1_cusp}
\end{center}

\newpage

\subsection{The case $c > 1$}
Now we get a new condition, see figures \ref{c4} and \ref{c4_cusp}:
\begin{equation}
	( \; A \le 0 \enspace \; and \enspace \;  B \le 0 \; )  \quad or \quad
   ( \; A \ge 0  \enspace \; and \enspace \; B \ge 0   \; \enspace \; and \enspace \; E \ge 0 \; )
   \label{Cond_cgt1}
\end{equation}

\begin{center}
  \includegraphics[width=0.92\textwidth]{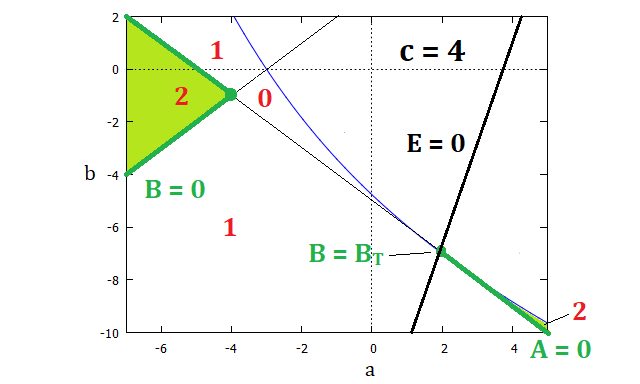} 
  \captionof{figure}{$c = 4$, roots of the cubic polynomial $x^3 + a x^2 + b x + c$,
                     in red the number of roots in the interval $[-1, 1]$.
                     The black line $E$ allows to distinguish the $2$ cases: $0$ roots in the interval
                     left of the line and the desired $2$ roots in the interval right of the line.
                     The line continues upwards to following figure with the cusp and passes there through
                     the point $A = A_T$.}
  \label{c4}
\end{center}

\begin{center}
  \includegraphics[width=0.92\textwidth]{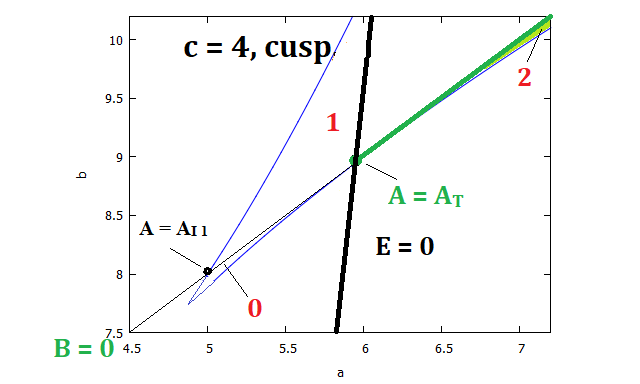} 
  \captionof{figure}{$c = 4$, the cusp, roots of the cubic polynomial $x^3 + a x^2 + b x + c$,
                     in red the number of roots in the interval $[-1, 1]$.
                     The black line $E$ allows to distinguish the $2$ cases: $0$ roots in the interval
                     left of the line and the desired $2$ roots in the interval right of the line.}
  \label{c4_cusp}
\end{center}

\newpage

\section{Cubic polynomials with $0, 1$ or $3$ roots in the interval}

To treat these remaining cases, the 5 quantities $A, B, A_T, B_T, E$ defined in  section \ref{inter} are sufficient.
The reader can find the conditions just looking at the figures and using other quadrants in the conditions.
The line $E$ is used in figure \ref{c4} in the right quadrant to distinguish between $0$ and $2$ roots in the interval. The same line $E$ can also be used in figure \ref{c1_4} in the upper quadrant to distinguish between $1$ and $3$ roots.

\section{Cubic polynomials with a pair of complex conjugated roots \\ and $0$ or $1$ roots in the interval}

These cases are simpler than the previous cases. Now the discriminant condition is $D_3 < 0$.
Only the $2$ quantities $A, B$ defined in  section \ref{inter} are sufficient for the conditions.
The case with different numbers of roots in a quadrant does not occur. \\

The previous figures contain enough information, to treat these cases too. Start from one of the $6$ regions with red numbers, $5$ of them in the following figure \ref{c1_4_cc}. For the $6^{th}$ region with a red $3$, see the detailed figure \ref{c1_4_cusp} with the cusp. Move towards the blue curve with $D_3 = 0$. For polynomials on this curve we get a real double root. Proceeding now into the region $D_3 < 0$ this real double root changes into a pair of complex conjugated roots. When $n$ is the number of roots in the interval in the starting region we get $n' = n \mod 2$ roots in $D_3 < 0$.

\begin{center}
  \includegraphics[width=0.92\textwidth]{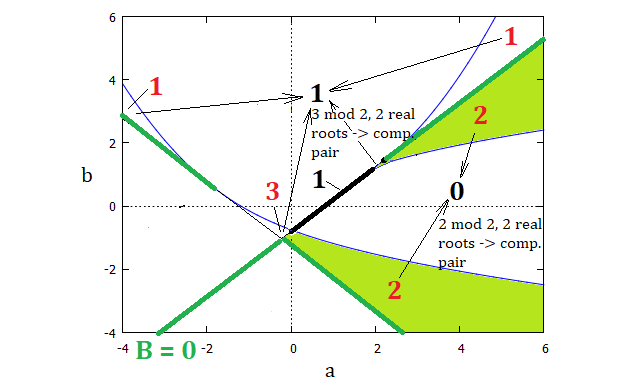} 
  \captionof{figure}{The here relevant parts of figure \ref{c1_4} for $c = 1 / 4.$
                     The black interval on the line $B = 0$ belongs to a polynomial with the
                     root $-1$ in $[-1, 1]$ and $2$ complex conjugated roots. The $4^{th}$ arrow
                     pointing to the black $1$ in the upper quadrant starts in the detailed figure
                     with the cusp.}
  \label{c1_4_cc}
\end{center}

The other cases for $c$ and the final result with the conditions are left for the reader.

\newpage

\section{Numerical and plausibility checks}

The results in section \ref{Cond_one_look} were checked numerically with thousands of cubic polynomials.
Rationals instead of floats for the coefficients were used. So it was possible, that the test also
covered polynomials on the dark green lines, belonging to the equalities $\ge, \le, =$ in the conditions. The test was also designed to cover the cases with double roots $D_3 = 0$ and with $c = 1$. \\

A plausibility check for the conditions: \\

The following $2$ maps generate a Kleinian 4-Group for $c \ne 0$.
The map $M : x  \rightarrow - x$ leaves 2 roots in $[-1, 1]$ in this interval
The map $N : x \rightarrow 1 /  x$ maps a root in  $[-1, 1]$ out of this interval, the other roots from outside into the interval. So N maps the problem "2 real roots in a closed interval" into the problem "1 root in an open interval". \\

The corresponding $2$ maps in the coefficient space are:
\begin{equation}
	M : a \rightarrow - a,   b \rightarrow b,   c \rightarrow - c  \qquad
	N : a \rightarrow b / c,   b \rightarrow a / c,    c  \rightarrow 1 / c
\end{equation}

Show how these maps N, M act on the $5$ quantities $A, B, A_T, B_T, E$, the discriminant $D_3$, and on the conditions \ref{Cond_c01}, \ref{Cond_cgt1}. Show how the maps change a condition from true to false.

\vspace{10mm}
\noindent \textbf{\large Appendices}

\appendix

\section{An example from physics: the Lagrange top  \label{ExPhysics}}

See \cite{HaFi}, chapter 3.6 "The Heavy Symmetric Top" (3.66) (3.72) and (3.73) with the cubic polynomial.

\begin{equation}
	a = \frac{I_3 \, \omega_3}{I_1}  \quad  b = \frac{p_\Phi}{I_1}  \quad  
   \alpha = \frac{2 \, E'}{I_1}  \quad  \beta = \frac{2 \, M \, g \, l}{I_1} 
   \label{top_abalbe}
\end{equation}
\begin{equation}
	(1 - u ^ 2) \, (\alpha - \beta \, u) - (b - a \, u) ^ 2 = 0
   \label{equ_nut}
\end{equation}
In the case of nutation the upper and lower limits for $\cos (\theta)$ are given by $2$ real roots of the cubic above in the closed interval $[-1, 1]$. The $3^{rd}$ real root lies outside this interval. With the conditions from section \ref{Cond_one_look} we get the following results. There should be no confusion with  the $a, b$ in \ref{top_abalbe} and the coefficients of the monic cubic depending on the context.

\begin{equation}
	A = - (a - b) ^ 2 / \beta  \quad B = - (a + b) ^ 2 / \beta  \quad
	A_T = - 4 \, (b ^ 2 - \alpha - \beta) / \beta \quad B_T = - 4 \, (b ^ 2 - \alpha + \beta) / \beta
   \label{top_AB}
\end{equation} \\

When $b \ne \pm a$ the $A, B$ in \ref{top_AB} are $A, B \ne 0$ and have the same sign. So the polynomial is located in the interior of the left or right quadrants (in the light green shaded region not on the
dark green lines on the boundary in the figures \ref{c0} and ff.). The coefficient $c$ of the monic
cubic is $c = (b ^ 2 - \alpha) / \beta$. \\

Let the discriminant condition $D_3 \ge 0$ be fulfilled. \\

Case 1: if $ \quad b \ne \pm a \quad $ and $ \quad -1 \le c \le +1 \quad $  then \quad $2$ roots in $[-1, 1]$. \\

There are $2$ remaining cases:
\begin{itemize}
   \item $b \ne \pm a$, $\vert c \vert > 1$, polynomials in the interior of the quadrants, the line $E$ is needed in a subcase
   \item $b = \pm a$, polynomials on the boundary of the quadrants with a root $\pm 1$, the intersection points $A_T, B_T$ are needed
\end{itemize}

They are left for the reader as exercise.

\bibliographystyle{amsplain}

\begin{thebibliography}{10}
\bibitem{Prod} E. M. Prodanov, \textit{On the cubic equation with its Siebeck-Marden-Northshield triangle and the quartic equation with its tetrahedron}. Journal of Computational Science 73 (2023)
 \href{https://doi.org/10.1016/j.jocs.2023.102123} {DOI: 10.1016/j.jocs.2023.102123}
\bibitem{HaFi} L. N. Hand and J. D. Finch, \textit{Analytical Mechanics}. Cambridge University Press, 1998
 \\ \href{https://www.damtp.cam.ac.uk/user/tong/dynamics/three.pdf} {Chapter 3. The Motion of Rigid Bodies}
\end{thebibliography}

\end{document}